\documentclass{article} 

\font\Bbb=msbm10
\def\R{\mbox{\Bbb R}} 
\def\D{{\cal D}}
\def\I{{\cal I}}
\def\O{{\cal O}}
\def\x{{\bf x}}
\def\rank{{\rm rank}}
\def\tomega{\widetilde{\omega}}

\usepackage{graphicx}
\usepackage{caption}

\def\wider{\advance \leftskip -\parindent \advance \rightskip -\parindent}
\def\plot#1#2#3#4#5{
        \begin{figure}[tb] %\begin{center}
	\setlength\captionmargin{-#5}
        \moveleft#5\hbox{\includegraphics[width=#3,height=#4,clip]{#1}}
        \caption{#2}
        \label{#1}
\end{figure}
}

\newcommand{\QED}{\hspace*{2em}\hfill$\bullet$}
\newtheorem{theorem}{Theorem}

\newtheorem{definition}[theorem]{Definition}

\def\smallmatrix{\null\,\vcenter\bgroup \baselineskip=7pt
  \ialign\bgroup\hfil$\scriptstyle{##}$\hfil&&$\;$\hfil
  $\scriptstyle{##}$\hfil\crcr}
\def\endsmallmatrix{\crcr\egroup\egroup\,}  

\begin{document}

\title{Area preservation in computational fluid dynamics}

\author{Robert I. McLachlan\\
\noindent $^{1}$ Mathematics, Institute of Fundamental
Sciences, \\
Massey University, Palmerston North, New Zealand; \\
email: R.McLachlan@massey.ac.nz}
        
\date{28 May 1999}
\maketitle

{
\begin{center}
\bf Abstract
\end{center}
\narrower\narrower

\noindent 
Incompressible two-dimensional flows such as the 
advection (Liouville) equation
and the Euler equations have a large family of conservation laws related
to conservation of area. 
We present two Eulerian numerical methods which preserve a discrete analog
of area. The first is a fully discrete model based on a rearrangement of
cells; the second is more conventional, but still preserves the area
within each contour of the vorticity field. Initial tests indicate that
both methods suppress the formation of spurious oscillations
in the field.

}

\section{Introduction}
When a smooth field $\omega(x,y)$ is advected by an area-preserving
flow, the area within each contour of $\omega$ is preserved. This is
seen in pure advection and in the Euler equations, for example, and is
important in the numerical solution of two-phase free boundary problems,
where the total volume of each fluid should be preserved.
Yet, although the advection problem has been addressed in probably 
thousands of papers, and very accurate, stable, and efficient methods
are known, no existing numerical methods take the area-preservation 
property into account.
In this Letter we present an initial study containing two methods
which do preserve (a discrete analog of) area. Although they are not,
presumably, competitive with the best existing methods for advection,
the results are extremely promising. 

The configuration space of an inviscid incompressible fluid is the
group $\D_\mu$ of volume-preserving diffeomorphisms of the fluid's
domain; the `Arnold' picture, in which the Euler equations
are geodesic equations on this group equipped with the kinetic
energy ($L_2$) metric, is treated in \cite{ar-kh}. 
The configuration at any time is a volume-preserving rearrangement
of the initial condition. 
Existing Eulerian numerical methods do not 
preserve any discrete analogue of this property. 
This is particularly relevant in two dimensions, where area
preservation leads to an infinite number of conserved quantities, 
the generalized enstrophies. 

We consider a two-dimensional fluid with divergence-free
velocity field ${\bf u}=(u,v)$, stream function $\psi$ (i.e.
$u = \psi_y$, $v=-\psi_x$), and some quantity $\omega$, which
we call the vorticity, which is advected by the fluid:
\begin{equation}
\label{advection}
\dot \omega + {\bf u}\cdot\nabla\omega = \dot\omega + J(\omega,\psi)  = 0,
\end{equation}
where the Jacobian
\begin{equation}
\label{jacobian}
J(a,b)={\partial(a,b)\over\partial(x,y)}.
\end{equation}
In the two situations we shall consider, this is a Hamiltonian
system with Poisson bracket
\begin{equation}
\label{pb}
\{F,G\} = \int \omega J(\delta F,\delta G)\, dx dy 
\end{equation}
and Hamiltonian
$$ H = {1\over2} \int \psi\omega \, dx dy $$
where the stream function $\psi$ is either a given function
$\psi = \psi(x,y,t)$, in which case Eq. (\ref{advection}) 
is the Liouville (advection)
equation,  or is
determined by the Poisson equation $\nabla^2 \psi = -\omega$,
in which case (\ref{advection}) is the 
the two-dimensional Euler equation. Other two-dimensional
flows such as the shallow water and semi-geostrophic equations also possess a
quantity $\omega$, called a potential vorticity, that is advected
according to Eq. (\ref{advection}). There are also applications to 
level-set methods, in which $\omega$ is not a physical variable but is
introduced so that the curve $\omega(x,y)=c$ can indicate a free boundary.

The Casimirs of the Poisson bracket (\ref{pb}) are conserved quantities
of the PDE (\ref{pb}). These can be variously written as
$$ C_f = \int f(\omega) \, dx dy$$
for any function $f$ such that $C_f$ exists, as
$$ C_n = \int \omega^n\, dx dy,$$
called the generalized enstrophies ($C_2$ is the usual enstrophy), or
as the areas enclosed by each vorticity contour
$$ A(c) = \int_{\omega\ge c}1\, dx dy.$$
They all reflect the fact that $\omega$ is being advected by
an area-preserving vector field and can only reach states
which are area-preserving rearrangements of its initial state.
That is,
$$ \omega(x,t) = \omega(\varphi_t^{-1}(x),0),$$
where $\varphi_t$ is the time-$t$ flow of the vector field ${\bf u}$. 

The famous Arakawa Jacobian is an Eulerian finite difference
approximation of Eq. (\ref{jacobian}) which preserves discrete analogues of
the energy $H$ and the enstrophy $C_2$ \cite{arakawa,mclachlan}. 
It is known to preserve the mean of the energy spectrum and to
prevent some nonlinear instabilities. However, the other conserved
quantities are not preserved and their role in the dynamics
is not known \cite{ar-kh}. 

Area preservation can also be studied in 
a Lagrangian framework---for example, point vortex methods could
be said to be area-preserving---but Lagrangian schemes carry a lot
of extra information (the particle paths) which should be
decoupled from the dynamics. The dimension of the phase space is
halved in Eulerian form, and further reduced by preserving
(discrete analogues of) the Casimirs. For ODEs, it is well
established that the best long-time results are obtained by working
in the smallest possible phase space \cite{mc-qu}. 

The Hamiltonian picture has been described by Marsden and
Weinstein \cite{ma-we}.
The configuration space is the 
group $\D_\mu$. The Euler equations in Lagrangian form are a canonical
Hamiltonian system on $T^*\D_\mu$, and in Eulerian form are
a Lie-Poisson system on the dual of the Lie algebra of $\D_\mu$,
which is identified with the space of vorticities. 
The coadjoint orbits of this space are the level sets of the
Casimirs, each of which is a symplectic manifold. Discretizations
of the Eulerian form are not, in general, Hamiltonian systems,
nor do they have conserved quantities corresponding to the
Casimirs (although there is one interesting Hamiltonian discretization,
the sine-Euler equations \cite{zeitlin}). 

Therefore we forget about the Hamiltonian structure and study
the Casimirs---the area-preservation---and present two
models in which a discrete analogue of the areas $A(c)$ is
preserved. The first (Section 2), based on a literal rearrangement of
cells, is interesting in that it gives a fully-discrete, 
cellular-automata-like model of an incompressible fluid. It
does not preserve smoothness of the vorticity field (although
filamentation and turbulence mean that it can't usually
stay very smooth anyway). A smooth version (Section 3) is based 
on computing an approximation of $A(c)$ which is smooth as 
a function of $c$, and relabelling the vorticity field so
that $A(c)$ is constant in time. It is tested on the Liouville
equation and prevents the appearance of large spurious 
maxima and minima in the vorticity field during its evolution.

\section{The cell rearrangement model}
Both of the models presented here are projection schemes. The
vorticity is evolved by any sensible scheme for some
short time $t$ (e.g., 1--10 time steps), and then projected onto
some space of rearrangements of the original vorticity.

In this section we consider the vorticity field to be piecewise
constant on a set of fixed cells, which for convenience we take
to be squares with side $h$. A (minuscule!) subset of the rearrangements of the
initial condition is given by the permutations of the cells. However,
these can be naturally associated with the fluid flow. For, consider
the area-preserving map $\varphi$ which is the time-$t$ flow of the
fluid. According to a theorem of Lax \cite{lax}, there is a mapping
$P$ which permutes cells and which satisfies 
\begin{enumerate}
\item $P(C)\cap\varphi(C)\ne\emptyset$ for all cells $C$; and
\item $\|P(x)-\varphi(x)\|\le \sup_{y,z\in C}\|\varphi(y)-\varphi(z)\|+
\sqrt{2}h\quad\forall x\in C.$
\end{enumerate}
The dynamics of such {\it lattice maps} are often studied. For example,
if the continuous map $\varphi$ is iterated on a computer, it will not be
exactly a bijection or exactly area-preserving, due to round-off error. By
replacing it with a lattice map and examining the limit $h\to0$ the
effects of roundoff error can be studied. 

The easiest way to construct lattice maps is as a composition of shears
$x_i' = x_i$ for $i=1,\dots,d$, $x_j' = \lfloor f_j(x_1,\dots,x_d)
\rfloor$ for $j=d+1,\dots,n$, where $\lfloor x\rfloor$ is 
the nearest lattice point to $x$. This would be suitable, for example,
if $\varphi$ itself were approximated by a product of shears,
as is for example the flow of separable Hamiltonians $H=H_1(p)+H_2(q)$
(the flow of the Hamiltonian vector field corresponding to each $H_i$
is a shear). This is very fast and the permutation need not be constructed
explicitly.

However, in the present case $\varphi$ can only be obtained by integrating
the Lagrangian particle paths for a short time $t$, and an explicit
approximating lattice map seems to be unobtainable. Scovel \cite{scovel}
suggested using maps of the form, e.g., $x' = x + \lfloor
J \nabla S((x+x')/2)\rfloor$ for a suitable Poincar\'e generating
function $S$ (here $S=(\Delta t)\psi$ would give a good approximation of
the time-$\Delta t$ flow of the stream function $\psi$). However, 
this nonlinear, discrete equation does not seem to have solutions in general.

Thus, it seems that one must laboriously construct a table of the
permutation. An algorithm which does this is described in \cite{mu-kl}.
Its running time is $\O(N^3)$, where $N=\O(1/h^2)$ is the number
of cells. One must construct lists of candidate cells (e.g.,
all those that intersect $\varphi(C)$) and make successive choices
from these lists, backtracking when no choices remain. While practical
for moderate $N$ when the dynamics of the lattice map are
going to be studied intensively, in the present application $\varphi$
changes at every time step; searching for a completely new permutation
every time is too expensive. This approach has been explored
by Turner \cite{turner}.

Luckily, there is a way out of this impasse, using the extra physical
information attached to each cell: the vorticity itself. The only
use of the permutation $P$ is to update the vorticity field $\omega$
by $\omega\mapsto\tomega$, $\tomega\circ P = \omega$,
in order that the distribution of vorticity values remains constant.
This can be achieved directly, without actually constructing
a $P$ which approximates $\varphi$, by the following algorithm.
Let $\rank_t(c)$ be the number of cells with vorticities
greater than or equal to $\omega$ at time $t$, i.e., 
$$ \rank_t(c)=\#\{j: \omega(x_j,t)\ge c\},$$
with ties broken arbitrarily to make $\rank_t$ an invertible function
onto $\{1,\dots,N\}$. Then:
\begin{enumerate}
\item Update the field $\omega$ for time $t$ any standard
Eulerian method; and
\item let $\tomega_j = \rank_0^{-1}(\rank_t(\omega_j))$.
\end{enumerate}
The new field $\tomega$ can be constructed
in time $N\log N$ by sorting the two lists of vorticity values
at times $0$ and $t$. The largest current value is replaced by
the largest original value, and so on. 
(Other updates, based on minimizing $\|\tomega-\omega\|$,
are also possible.)

This algorithm
can be regarded as constructing a permutation, albeit a permutation
that has no relationship to the flow $\varphi$.
This is a truly finite-state model of an incompressible fluid: the
state space is the permutation group $S_N$ and the fluid dynamics
reduces to the dynamics of the map $S_N\to S_N$, defined above,
parameterized by the initial distribution of vorticity values.
It is an almost cellular-automata-like fluid model, although lacking
the local update property of CAs. It has the aesthetic appeal of capturing
the vorticity-rearrangement property perfectly in a naturally 
discrete way, of constructing a ``discrete coadjoint orbit'',
and it is very cheap.

However, these advantages are
offset by a practical disadvantage of lack of smoothness.
The new vorticity values are selected somewhat arbitrarily from
the sorted list, and the new field may be rougher
than the original. This is probably unavoidable, given the chosen
fully discrete state space. The noise of this imposed roughness
may swamp any gains from preserving the coadjoint orbits.
However, in a turbulent flow with highly filamented vorticity,
the loss of smoothness may not be significant. A second
consequence of the lack of smoothness is that if $|\omega(x,t)
-\omega(x,0)|$ is too small, then $\tomega\equiv\omega$---the
field cannot be updated at all. The remapping interval $t$ must
be large enough to allow some change in the configuration. 
For example, the flow map $\varphi$ should move each cell across at
least 2 cells so that the algorithm has some scope for finding
a suitable permutation. 

\section{The vorticity relabelling model}
The cell rearrangement model produces an area function $A(c)$ which
is discon\-tin\-uous---in fact, it is piecewise constant. To improve
it, we need to
\begin{itemize}
\item[(i)]
produce a smoother approximation of $A(c)$. If
$\omega(x,y)$ is a smooth function, we want an approximation
of $A(c)$ which is as smooth and accurate as possible, using
only the grid values $w(x_i,y_j)$; and
\item[(ii)]
project the vorticity function so that its area function
$A(c)$ at time $t>0$ equals (or closely approximates) the
initial area function.
\end{itemize}

\subsection{Computing the areas enclosed by vorticity contours}
We consider a compact domain $\Omega$ with area 1, 
usually a square or torus, on which
$\omega$ is bounded with range 
$[\omega_{\rm min},\omega_{\rm max}]$, and of smoothness $C^r$. 
It may be degenerate, e.g., constant on open sets.
\begin{definition}
The {\em area function} of the field $\omega$ is 
the area enclosed by the set $\{(x,y): \omega(x,y)\ge c\}$, i.e.,
$$ A_{\omega}:[\omega_{\rm min},\omega_{\rm max}]\to[0,1],\quad
A_{\omega}(c) = \int_{\omega(x)\ge c}1\, dx\,dy $$
\end{definition}
$A(c)$ is strictly decreasing with respect to $c$.
It is $C^r$ at regular (noncritical) values $c$,
$C^0$ at nondegenerate critical values, and discontinuous
at $c$ if the set $\{x:\omega(x,t)=c\}$ has positive area.
(Lack of differentiability at critical values can be seen by studying
$\omega = -(x^2 + y^2)$, for which 
$A(c) = \pi c$ for $c\le 0$ and $0$ for $c>0$.) 
Thus, its inverse $A^{-1}$ exists and is nonincreasing
(i.e., more area must be enclosed by a lesser value of $\omega$.)

Let $\I$ be an interpolation or approximation operator
mapping grid functions to fields, i.e. functions on $\Omega$. 
\begin{definition}
The area function of the grid function $\omega$ is 
defined to be the area function of its interpolant, i.e.,
$$ A_{\omega} := A_{\I\omega}. $$
\end{definition}
It automatically inherits the monotonicity
properties of $A$. Let ${\cal R}$ be a restriction operator mapping
fields to grid functions, usually by evaluating on the grid.
Let $\tomega = \I{\cal R}\omega$.
The crucial observations are the following:
\begin{enumerate}
\item Choice of $\I$ can lead to $A_{\tomega}$ being as
smooth as $A_{\omega}$, and of any order of accuracy as
an approximation;
\item If $\tomega$ is piecewise linear, its contours
are polygons, whose area can be found quickly for any contour topology;
\item If $\tomega$ has polygonal contours, $A_{\tomega}$
can be second order accurate and as smooth as $A_{\omega}$.
\end{enumerate}
Item (1) is obvious, and is a consequence of existence of 
$C^r$ approximations to functions. An algorithm for finding
the area enclosed by (unions of) polygons is given below.
The most important point is (3), as it says that smooth
interpolants, which are expensive in two dimensions, are not
needed to compute a smooth area function.

Consider a grid function on a triangulation of $\Omega$. Interpolating
by piecewise polynomials along edges only, and constructing
the interpolant whose contours are line segments within each
triangle (whose graph is a ``ruled surface''), 
yields an area function which is 
as smooth as the interpolant at vertex (grid point) values and analytic
elsewhere. Thus, only smooth one-dimensional interpolation
is needed, which is relatively cheap (e.g., $C^1$ can be
achieved using local cubics).

Piecewise linear interpolation yields a $C^0$, second-order-accurate
area function. However, it is better than its mere continuity might
make it appear, since
its derivative jumps at vertex values are only ${\cal O}(h^2)$
on a grid with spacing $h$. So,
numerically, it is indistinguishable from a $C^1$ function.
In practice, the most glaring jumps are in its {\it second}
derivative at vertex values, not in the function itself.
(See Fig. 2.)
Piecewise linear interpolation seems to be suitable in practice and
this is what we use in the tests below.
(If the main computational
grid is square, we triangulate using an extra vertex at the
center of each cell, whose value is assigned by linear interpolation.)
However, true $C^1$ area functions have been tested as well.

The great advantage of polygonal contours is that
the area of a simple polygon with vertices $\x_i$, $i=1,\dots,n$, is
very easy to compute: it is
$${1\over2}\sum_{i=1}^n \x_i\times \x_{i+1},$$
where $x_{n+1}:=x_1$. This can be seen by deriving it for a triangle,
triangulating the polygon against a fixed point, and then using
independence with respect to the fixed point. It can also be
viewed as a discretization of
$$ \int_\Omega 1\, dx\, dy = {1\over2}\int_\Omega d(x dy - y dx)
= {1\over2} \int_{\partial\Omega}\x \times d\x.$$
However, it would be expensive to chase contours around the
domain and construct a list of simple polygons. Instead, one
can simply scan each triangle for occurrence of a contour,
find its endpoints $\x_1$, $\x_2$, and accumulate
$\x_1\times\x_2$ with a sign determined by the sense of the 
triangle when its vertices are visited in order of
increasing function values. This handles arbitrary contour topology.
(Exception handling is needed when two vertices and the contour
all have equal values.)

We are not sure if there is a similarly simple method with
higher order contours.  In practice, to get more than
second order accuracy, we use Richardson extrapolation from
a coarser grid.

For a list of contour values, 
the above algorithm involves scanning the cells once and accumulating
areas of
the relevant contours. If $N_c$ values of the area function
are needed, and the grid size is $\O(h)$, then each cell will 
contain $O(h N_c)$ contours on average,
so the computation takes time $\O(N_c/h)$. In practice, we
take $N_c=\O(1/h)$ and build a function table, which is
later interpolated as needed. Thus computing the areas
takes $\O(1/h^2)$, i.e., it is linear in the number of grid points.

If $\omega$ is nearly constant on large areas, then $A(c)$ can
be very steep, so it should be tabulated using adaptive stepping
in $c$. An example of the $C^0$ estimate of $A(c)$ given by piecewise
linear interpolation is shown in Fig. 2, together with the piecewise
constant estimate given by simply counting the number of vertices
where $\omega>c$. Its (numerical) derivative indicates its smoothness.
Some care must be taken when interpolating to maintain monotonicity.

We have also computed smooth approximations of $A_\omega(c)$ for 
random $\omega$ fields whose contours have complicated topology.

\subsection{Projecting the the space of rearrangements}

After evolving $\omega$ for a short time $t$ with Eulerian method, 
we have two grid functions,
the vorticity at time $0$, $\omega(0)$, and at time $t$, $\omega(t)$. 
We wish to project $\omega(t)$ so that it is an (approximation of)
a rearrangement of $\omega(0)$. The projection should be small
and should not destroy smoothness. Traditional methods for
enforcing constraints, such as steepest descents, appear to
be completely infeasible because of the global and sensitive
dependence of $A(c)$ on the vertex values of $\omega$. Our
proposed method is a continuous version of the sorted-assignment
used in the cell rearrangement model of Section 2. In words, we compute the
area enclosed by the contour through each vertex value and
replace it by the value that originally enclosed that much area.
The contour shapes and topologies do not change: only
the values associated with each contour change.

\begin{definition}
The relabelling projection on grid functions $\omega_i\approx
\omega(x_i,t)$ is
defined by $\omega(t)\mapsto \tomega(t)$, where
\begin{equation}
\label{relabel}
A_{\omega(0)}(\tomega_i) = A_{\omega(t)}(\omega_i)
\end{equation}
for each vertex $i$.
\end{definition}
It is well defined by monotonicity of $A(c)$. It has an obvious
continuum analog (replacing $i$ by $x$ in Eq. (\ref{relabel})), 
which if applied to
every value of $\omega$ taken by a smooth vorticity field, 
with $\omega(0)$ and $\omega(t)$ both $C^r$, yields a new field
$\tomega$ that is $C^r$ away from critical points of $\omega_0$ and 
$\omega_t$ and $C^0$ at such critical points.

To compute a good approximation of this projection quickly, the current area
function is tabulated and interpolated at the vertex values,
and then $A^{-1}_{\omega(0)}$ (which, of course, does not change
during the run) is evaluated by interpolation. Of course,
we do not have a true projection in that $\widetilde{\tomega}
\ne\tomega$, because interpolation errors in the contours
do change the contour shapes by a small amount when the
vertex values are changed. We do not quite get
$A_{\omega(0)}(c)= A_{\tomega(t)}(c)$ for all $c$.
However, these errors can be
controlled independently of the discretization error in
$\omega$, for example, by using a higher order approximation
of $A$. In a numerical test, one application of Richardson 
extrapolation to the areas enclosed by piecewise linear 
contours gave $|A_{\omega(0)}-A_{\tomega(t)}|\sim 10^{-4}$
on a relatively coarse $20\times20$ grid. By contrast, without 
the relabelling projection, errors in the area function
rapidly reach order 1.

If the vorticity is evolved for a short time $t$, with a method
of spatial order $p$, spatial errors dominate the error in the 
area function which are 
$\O(t h^p)$. Thus, with $t=o(1)$, the projection only alters the
field by $o(h^p)$, and the overall method (after
evolution and projection), is still consistent
of the same order $p$. The projection cannot correct any errors
in the shapes of the contours, but it can stop those errors
growing further by propagation of the false distribution
of vorticity values, which is particularly bad for the 2D Euler
equations, where those values determine the velocity field itself.

\section{Numerical tests}
Here we illustrate some short tests to validate our approach and show
that it is indeed possible to compute and preserve area in an 
Eulerian method. We use a coarse ($20\times 20$) grid which barely
resolves the solution, and a crude (second order) finite difference
approximation to the spatial differences, in order to test whether the 
method can correct the large oscillations and area errors that result.

We solve the Liouville equation in $\Omega=[0,1]^2$
with initial field
$\omega = \exp(-45(x-{3\over4})^2 -15(y-{1\over2})^2)$ 
advected by the velocity field with stream
function $\psi = \sin(\pi x)\sin(\pi y)$. (See Figure 1.)
This velocity field has shear, so $\omega$
rapidly rolls up into a tight spiral, mimicking 
the filamentation of vorticity in the Euler equations. The spatial
derivatives in Eq. (\ref{advection}) are approximated by the
Arakawa Jacobian, which is second order and preserves discrete analogues
of energy and enstrophy. Although the discrete enstrophy $\sum\omega_i^2$
is preserved, this does not help the scheme preserve areas any better
than (nonconservative) central differences do.

Particles at the maximum of $\omega$ have a period of about 0.75. We integrate
with a second order method for 400 times steps of $\Delta t = 0.003$, or
total time $1.2$, during which this maximum rotates 1.6 times around the
centre of the square $\Omega$. 
Spatial errors completely dominate the total error at $t=1.2$.

Without any projection, oscillations rapidly develop and the distribution
of vorticity values is not maintained at all well (see Figure 3(a)). 
A large minimum of $\psi = -0.69$ forms, next to a spurious local maximum 
of $\psi = 0.46$.  The initial maximum of 1 has not been preserved but has 
decayed to 0.87. The comparison between the initial and final area functions
(see Figure 2) shows that the area within most vorticity contours is not
preserved at all.

The area-preserving methods both involve periodically remapping the vorticity.
If this period is too short (e.g. one time step), then the cell rearrangement
model cannot update the vorticity at all. If it is too long, then not only
the area but also the topology of the level sets can alter, which can
not be corrected by the present methods. Once a small island of vorticity
has been created, for example, it must be advected by the flow. 

We first consider the cell rearrangement model of section 2. Suppose
the remapping is applied
every $N_r$ time steps. This needs a large $N_r$ to yield a reasonably
smooth remapped vorticity field; but if $N_r$ is too large then spurious
maxima can evolve which are not removed by the remapping. With no remapping,
this maximum reaches $0.46$. With $N_r=50$, it reaches $0.26$. With
$N_r=20$, there is no isolated spurious maximum, but oscillations start
to appear within the main island. These grow worse at $N_r=10$. Therefore,
$N_r=20$ seems a reasonable balance, and the final field is shown in 
Figure 3(b). In one remapping period, the central peak moves across
about 2 cells.

This remapping is very fast, but it does not maintain smoothness of $\omega$,
as can be seen here. In fact, it is surprising that it works even as well
as it does in this example. However, the lack of smoothness would not 
be a problem in problems involving poorly-resolved turbulent fields. 

We consider now the vorticity relabelling model of section 3.
In this model we are free to decrease the remapping interval $N_r$
as desired: we still obtain smooth results with $N_r=1$, for example.
As $N_r$ is decreased, the results progressively improve.
For $N_r=\infty$, 20, 10, and 5, the peak of the spurious maximum is at
$\omega = 0.46$, $0.09$, $0.02$, and $0.002$, respectively. (Because
of its smooth interpolation, it cannot completely eliminate this maximum,
as the cell rearrangement model does.) Results for $N_r=10$ are shown
in Figure 3(d). The final field is very smooth, considering the coarse
$20\times 20$ grid, and very plausibly represents an element of 
the original state composed with an area-preserving diffeomorphism.
One contour of the exact solution (found by particle tracking) is shown
in the background. The computed solution has clearly suffered far too
much diffusion, a result of using diffusive, non-upwinded second
differences to approximate the advection term. Nevertheless, it is
impressive that such information can be extracted from the same
method that produced Figure 3(a), by merely imposing some conservation
laws. 

Finally, Figure 3(c), shows the vorticity relabelling model
applied to an even simpler spatial discretization, namely ordinary central
differences. It is in fact {\it more} accurate that the Arakawa Jacobian
(Fig. 3(d)), being slightly less diffusive. Thus, preserving areas lets one
use much simpler finite differences and still maintain smooth, non-oscillatory
solutions.

Any of the techniques presented here can be combined with a more
sophisticated underlying Eulerian scheme. If we used a high-order,
low-diffusion upwinding scheme, for example, then area errors would have
been much less than in Fig. 3(a); but they would still increase over time.
Applying the vorticity relabelling would still improve the solution.

\section{Discussion}

The methods discussed here take into account one large family of conservation
laws. This possibility raises many questions. What is the effect of using
these methods for very long times?
What is their effect on other conservation laws
such as energy and symplecticity? 
How well do they work on larger applications such as the shallow water
equations? (For level-set applications, a simpler update, adding
a constant to the advected field so that the area inside one particular
level set is preserved, may be preferable.)

More theoretically, is it possible to regard the `equal
area' functions as defining a discrete phase space in which consistent
approximations can be directly derived, instead of using brute force
modification of existing methods?
While desirable, this looks difficult,
since we are not projecting to any well-defined manifold. Consider
the subset of $\R^{N^2+1}$ defined by
$$ A_{\omega}(\omega_i)=A_0(c),\ i=1,\dots,N^2.$$
This does have dimension 1 in general, but is formidably
curled up on itself. It may be better to think of the constrained phase
space as the configurations lying within some small distance of 
a manifold of dimension $\R^{N^2-\O(N)}$, as we are enforcing one curve's
worth of constraints.

\bigskip
{\small
\noindent{\bf Acknowledgements}
I am extremely grateful to Tom Hou and Arieh Iserles for bringing this
problem to my attention, to Reinout Quispel for useful discussions
and for providing the reference \cite{mu-kl}, and to Paul Turner who
studied the cell rearrangement model in the course of his M.Sc. 
thesis \cite{turner}. This research is supported
by the Marsden Fund of the Royal Society of New Zealand. Part of it 
was undertaken when the author enjoyed the support of the MSRI, 
Berkeley. MSRI wishes to acknowledge the support of the NSF through grant 
no. DMS--9701755.

}

\plot{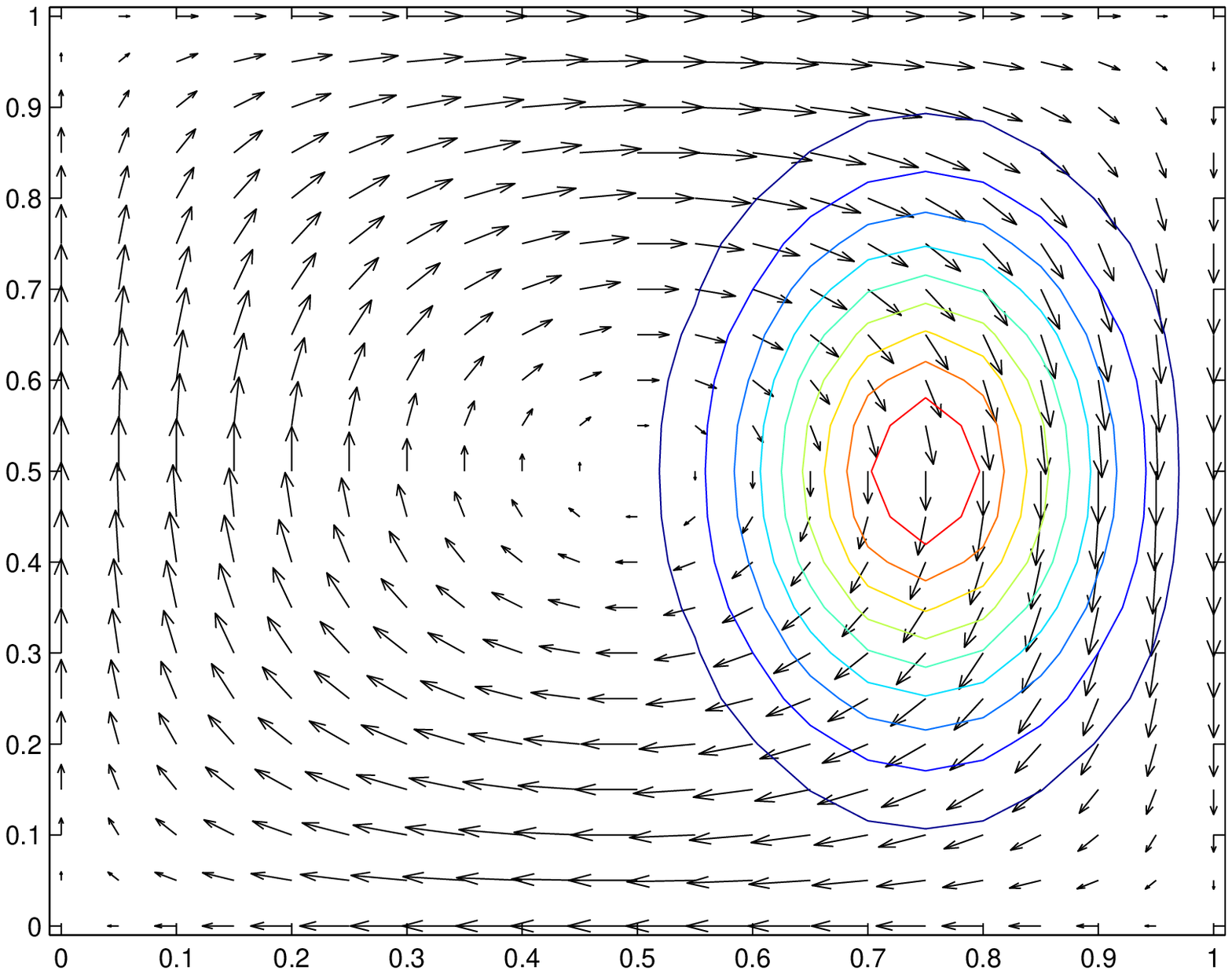}{Initial condition for the test problem in Section 4. The 
contours show level sets of $\omega(0) = \exp(-45(x-{3\over4})^2
- 15(y-{1\over2})^2)$. The arrows show the vector field corresponding
to the stream function $\psi = \sin(\pi x)\sin(\pi y)$ by which $\omega$
is advected.}{\hsize}{\hsize}{0cm}

\plot{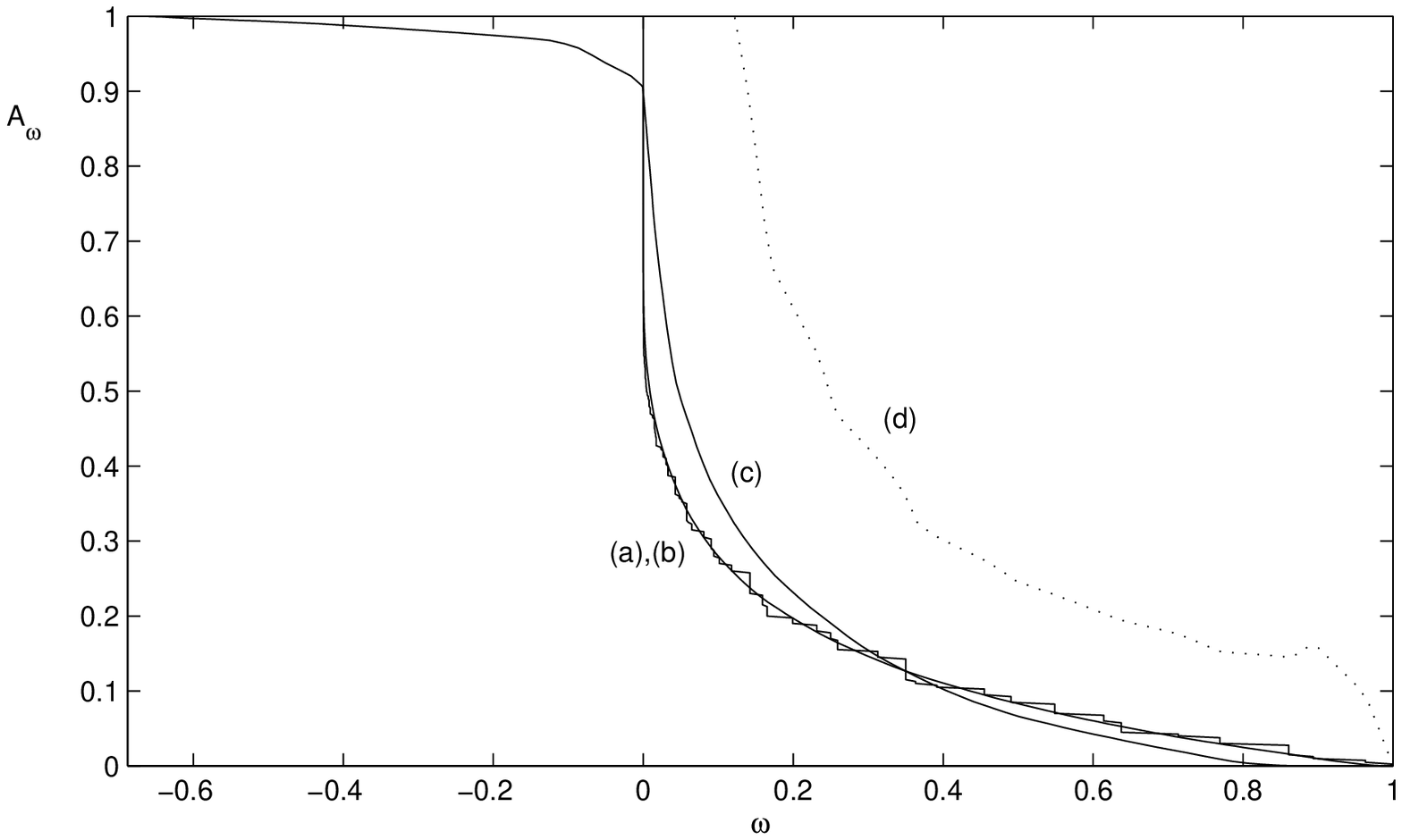}{Numerical computation of the area enclosed within vorticity
contours. (a): $C^0$ approximation to $A_{\omega(0)}$ using piecewise 
linear interpolation (Section 3). Here $\omega(0)$ is the initial condition
shown in Figure 1. (b): Piecewise constant approximation to $A_{\omega(0)}$ by 
sorting the list of vorticity values (Section 2). (c): Area function
$A_{\omega(t)}$
after evolving for time $t=1.2$ with no area preservation with an enstrophy-%
preserving scheme (Section 4). In the vorticity relabelling projection,
vorticity values are mapped from this curve back to (a). (d): Finite
difference approximation to $dA_{\omega(0)}(c)/dc$, showing that, although
only $C^0$, for numerical purposes it can be regarded as being differentiable.
The kinks in this derivative are due to $\omega$ being set
to zero on the boundary.}{16cm}{10cm}{1.94cm}

\plot{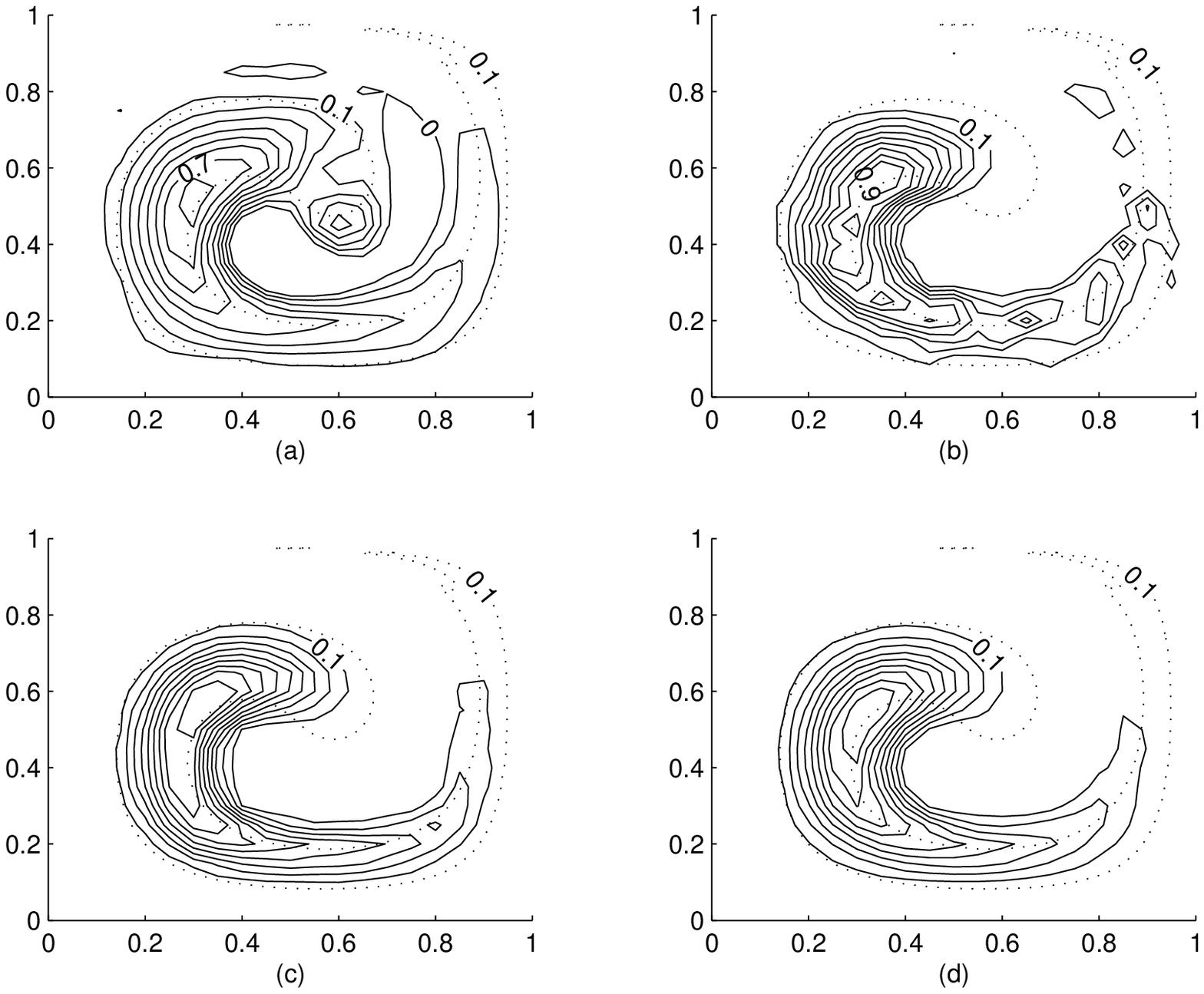}{Results for the advection problem on Fig. 1 after 1.6 rotations
about the center. (a): Arakawa differences with no area preservation. A
large negative blob of vorticity forms and spawns a secondary positive
blob. The dotted contour indicates the exact solution. (b): Arakawa
differences with cell rearrangement applied every 20 time steps ($\Delta t
= 0.001$). (c): Central differences with vorticity relabelling applied
every 10 time steps. (d) Arakawa differences with vorticity relabelling
applied every 10 time steps.}{16cm}{16cm}{1.94cm}

\end{document}